\documentclass[11pt,a4paper]{article}

\usepackage{epic,eepic,epsfig,amsmath,amssymb,amsthm}
\usepackage{t1enc}
\usepackage[francais,english]{babel}

\usepackage{caption}

\usepackage{color}

\usepackage{bbm}
\def\virgp{\raise 2pt\hbox{,}}

\renewcommand{\leq}{\leqslant}
\def\N{{\mathbb N}}

\def\R{{\mathbb R}}

\def\virgp{\raise 2pt\hbox{,}}
\def\cdotpv{\raise 2pt\hbox{;}}

\def\1{\mathbbm{1}}

\theoremstyle{remark}
\newtheorem{remark}{Remark}[section]

\theoremstyle{definition}

\theoremstyle{definition}

\theoremstyle{definition}

\begin{document}

\title{The Finite difference method for the Minkowski Curve}

\author{Nizare Riane, Claire David}

\maketitle
\centerline{Sorbonne Universit\'es, UPMC Univ Paris 06}

\centerline{CNRS, UMR 7598, Laboratoire Jacques-Louis Lions, 4, place Jussieu 75005, Paris, France}

\begin{abstract}
In this work, we describe how to approximate solutions of some partial differential equations using the finite difference method defined on the Minkowski self-similar curve. 
%We also give a numerical estimate of the error as a function of the graph order.
\end{abstract}

\maketitle
\vskip 1cm


\vskip 1cm

\vskip 1cm

\vskip 1cm

\section{Introduction}

\hskip 0.5cm The last few decades saw the birth of analysis on fractals, particularly, differential equations on fractals, introduced by J.~Kigami~\cite{Kigami1989},~\cite{Kigami1993},~\cite{KigamiStrichartzWalker2001},~\cite{Kigami2003}. May one consider a fractal set~$\cal F$, which shows enough nice properties (self-similar, post-critically finite), and a continuous function~$u$ defined on~$\cal F$, a Laplacian~$\Delta\,u$ can be obtained as the renormalized limit of a sequence of graph Laplacians~$\Delta_m u$,~$m\,\in\,\N$, built on a sequence of graphs which approximate~$\cal F$. \\

The seminal work of J.~Kigami has been followed by numerous contributions, in the field of fractal analysis. One may now explore the theoretical and numerical related areas with specific tools. But there still remain many points that ought to be studied in terms of numerical computation on specific fractals.\\

The door was opened by K.~Dalrymple, R.~S.~Strichartz, and J.~Vinson~\cite{DalrympleFDM}, who, in the case of he Sierpi\'{n}ski gasket~$\cal SG$, give an equivalent method for finite difference approximation. This work was followed by the one of M.~Gibbons, A.~Raj and R.~S.~Strichartz~\cite{GibbonsFEM}, where they describe how one can build approximate solutions, by means of piecewise harmonic, or biharmonic, splines, again in the case of~$\cal SG$.  They go so far as giving theoretical error estimates, through a comparison with experimental numerical data. One may note that this method is quite general, and can be extended to other self-similar sets.\\

In~\cite{RianeDavidM}, we built a Laplacian on the Minkowski curve~$\mathfrak MC$, which had seemed of interest for several reasons. First, it was a curve, and the challenge, as regards fractal ones, is to consider energy forms more sophisticated than classical ones, by means of normalization constants that could, not only bear the topology, but, also, the geometric characteristics.  One must of course bear in mind that a fractal curve is topologically equivalent to a line segment, and that Dirichlet forms, required for who aims at building a Laplacian, solely depend on the topology of the domain, and not of its geometry. Second, we came across applications in electromagnetism, in the field of fractal antennas, where one requires accurate numerical results.\\

%The Minkowski curve $\mathfrak MC$is a self-similar curve with a strange geometry, it has a fractal dimension of~$\displaystyle \frac{3}{2}$.

Our Laplacian on~$\mathfrak MC$ is obtained explicitly, as a limit of difference quotients. This opens the way to a theoretical and numerical study of partial differential equations, as the heat or wave equation. We hereafter present the resulting finite difference scheme. We did not find any equivalent in the existing literature.  \\

% $$\Delta u(x)=\lim_{m\rightarrow\infty} 64^m \left( \sum_{x \underset{m}\sim y}  u(kh,y)- u(kh,x)\right)$$

We proceed as follows: in section $1$, based on our previous results~\cite{RianeDavidM}, we present a mathematical construction of the Minkowski curve. In section~$2$, we successively build the finite difference schemes for the heat and the wave equations ; we then give a numerical estimate for the error. Computations have been made using Mathematica.\\

\section{The Minkowski Curve}

\noindent In the sequel, we recall results that are developed in~\cite{RianeDavidM}.\\

\noindent We place ourselves in the Euclidean plane of dimension~$2$, referred to a direct orthonormal frame. The usual Cartesian coordinates are~$(x,y)$.\\

\noindent  Let us denote by~$P_0$ and~$P_1$ the points:

$$P_0= (0,0) \quad , \quad P_1= (1,0)$$

\vskip 1cm

\noindent Let us denote by~$\theta \,\in\,]0,2\,\pi[$, $k>0$,~$T_1$, and~$T_2$ real numbers. We set:
$${\cal R}_{O,\theta}=\left(
\begin{matrix}
  \cos  \theta  & -   \sin  \theta \\
 \sin  \theta  &   \cos  \theta \\
\end{matrix}
\right)$$

\noindent We introduce the iterated function system of the family of~maps from~$\R^2$ to~$\R^2 $:
$$\left \lbrace f_{1},...,f_{8} \right \rbrace$$
\noindent where, for any integer~$i$ belonging to~\mbox{$\left \lbrace 1,...,8 \right \rbrace$}, and any~$X \,\in\R^2$:
$$f_i(X)=k\, {\cal R}_{O,\theta }
X +
\left(
\begin{matrix}
T_1\\
T_2\\
\end{matrix}
\right)
$$

\vskip 1cm

\begin{remark}
\noindent If $0<k<1$, the family~$\left \lbrace f_{1},...,f_{8} \right \rbrace$ is a family of contractions from~$\R^2$ to~$\R^2$, the ratio of which is~$k$.

\end{remark}

\vskip 1cm

\noindent According to~\cite{Hutchinson1981}, there exists a unique subset $\mathfrak{MC} \subset \R^2$ such that:
\[\mathfrak{MC} = \underset{  i=1}{\overset{8}{\bigcup}}\, f_i(\mathfrak{MC})\]
\noindent which will be called the Minkowski Curve.\\
\noindent For the sake of simplicity, let us set:

$${\cal F}=\displaystyle \underset{  i=1}{\overset{8}{\bigcup}}\, f_i  $$

\vskip 1cm

\noindent We will denote by~$V_0$ the ordered set, of the points:

$$\left \lbrace P_{0},P_{1}\right \rbrace$$

\noindent The set of points~$V_0$, where, for any~$i$ of~\mbox{$\left \lbrace  0,1  \right \rbrace$}, the point~$P_1$ is linked to the point~$P_{2}$, constitutes an oriented graph, that we will denote by~$ {\mathfrak {MC}}_0$.~$V_0$ is called the set of vertices of the graph~$ {\mathfrak {MC}}_0$.\\

\noindent For any strictly positive integer~$m$, we set:
$$V_m =\underset{  i=1}{\overset{8}{\bigcup}}\, f_i \left (V_{m-1}\right )$$

\noindent The set of points~$V_m$, where the points of an $m$-cell are linked in the same way as ${\mathfrak {MC}}_0$, is an oriented graph, which we will denote by~$ {\mathfrak {MC}}_m$.~$V_m$ is called the set of vertices of the graph~$ {\mathfrak {MC}}_m$. We will denote, in the following, by~${\cal N}_m$ the number of vertices of the graph~$ {\mathfrak {MC}}_m$.

\vskip 1cm

\noindent For numerical purpose, we use the following similarities, from~$\R^2$ to~$\R^2$, such that, for any~$X\,\in\,\R^2$:

$$
f_1(X) =\displaystyle \frac{1}{4} \,\left( {\cal R}_{O,0 }\,X +
\left(
\begin{matrix}
0\\
0\\
\end{matrix}
\right)
\right) \quad , \quad
f_2(X) =\displaystyle\frac{1}{4}\, \left( {\cal R}_{O,\frac{\pi}{2}}\,X +
\left(
\begin{matrix}
1\\
0\\
\end{matrix}
\right)
\right)  \quad , \quad
f_3(X) =\displaystyle\frac{1}{4} \,\left({\cal R}_{O,\frac{\pi}{2}}\,X +
\left(
\begin{matrix}
1\\
1\\
\end{matrix}
\right)
\right) $$

$$
f_4(X) =\displaystyle\frac{1}{4} \,\left({\cal R}_{O,\frac{3\,\pi}{2}}\,X +
\left(
\begin{matrix}
2\\
1\\
\end{matrix}
\right)
\right)  \quad , \quad
f_5(X) =\displaystyle\frac{1}{4} \left({\cal R}_{O,\frac{3\,\pi}{2}}\,X +
\left(
\begin{matrix}
2\\
0\\
\end{matrix}
\right)
\right)  \quad , \quad
f_6(X) =\displaystyle\frac{1}{4} \,\left({\cal R}_{O,0}\,X +
\left(
\begin{matrix}
2\\
-1\\
\end{matrix}
\right)
\right) $$

$$
f_7(X) =\displaystyle\frac{1}{4}\, \left({\cal R}_{O,\frac{\pi}{2}}\,X +
\left(
\begin{matrix}
3\\
-1\\
\end{matrix}
\right)
\right)  \quad , \quad
f_8(X) =\displaystyle\frac{1}{4}\, \left({\cal R}_{O,0}\,X +
\left(
\begin{matrix}
3\\
0\\
\end{matrix}
\right)
\right) \\
$$

\vskip 1cm

\section{The finite difference method on fractal sets}

In the sequel, we will denote by~$T$ a strictly positive real number, by~$M$ the number of contractions, by~$N_0$ the cardinal of $V_0$, and by~${\cal F}$ the considered fractal set.\\

\subsection{The heat equation}

\subsubsection{Formulation of the problem}

\noindent We may now consider a solution~$u$ of the problem:

$$ \left \lbrace  \begin{array}{ccccc}
\displaystyle \frac{\partial u}{\partial t}(t,x)-\Delta u(t,x)&=&0 & \forall (t,x)\,\in \,\left]0,T\right[  \times {\cal F}\\
u(t,x)&=&0 & \forall \, ( x,t) \,\in \,\partial {\cal F} \times \left[0,T\right[\\
u(0,x)&=&g(x)  & \forall \,x\,\in \,{\cal F}
\end{array}\right.$$

\noindent In order to define a numerical scheme, one may use a first order forward difference scheme to approximate the time derivative~$\displaystyle \frac{\partial u}{\partial t}$. The Laplacian is approximated by means of the graph Laplacians~$\Delta_m \,u$, defined on the sequence of graphs~$\left ( {\mathfrak MC}_m\right)_{m\in\N}$. \\

\noindent To this purpose, we fix a strictly positive integer~$N $, and set:

$$ \displaystyle{h =\frac{T}{N}}$$

\noindent One has, for any integer~$k$ belonging to~$\left \lbrace 0, \hdots, N-1 \right \rbrace$:

$$   \forall \,X\,\in \, {\cal F} \,:   \quad
\displaystyle\frac{\partial u}{\partial t}(kh,x)=\displaystyle\frac{1}{h}\,\left( u((k+1)h,X)-u(kh,X)\right)+{\cal O}(h)
$$

\noindent According to~\cite{RianeDavidM}, the Laplacian on the Minkowski Curve~$\mathfrak MC$ is given by:

 $$ \forall \, X\,\in \,{\cal F} \,:   \quad
\Delta u(t,x)= \lim_{m\rightarrow +\infty} 64^m \, \left( \sum_{X \underset{m}\sim Y}  u(t,Y)- u(t,X)\right)
$$

\noindent This enables one to approximate the Laplacian, at a~$m^{th}$ order,~$m \,\in\,\N^\star$, using the graph normalized Laplacian as follows:

$$ \forall \, k\, \left \lbrace 0, \hdots, N-1 \right \rbrace,\,\forall \,X\,\in \,{\cal F} \,: \quad
\Delta u(t,X)\approx 64^m \,\left( \sum_{X \underset{m}\sim Y}  u(k\,h,Y)- u(k\,h,X)\right)
$$

\noindent By combining those two relations, one gets the following scheme, for any integer~$k$ belonging to~$\left \lbrace 0, \hdots, N-1 \right \rbrace$, any point~$P_j$ of~$V_0$,~$0\leq j \leq N_0-1$, and any~$X$ in the set~\mbox{$V_m\setminus V_0$}:

$$\left ({\cal S}_{\cal H}\right) \quad \left \lbrace \begin{array}{cccc}
\displaystyle \frac{u((k+1)\,h,X)-u(k\,h,X)}{h}&=&64^m \,\left(\displaystyle \sum_{X \underset{m}\sim Y}  u(k\,h,Y)- u(k\,h,X)\right) &   \\
u(k\,h,P_j)&=&0 &     \\
u(X,0)&=&g(X) &
\end{array} \right.$$

\noindent Let us define the approximate equation as:

$$
u((k+1)\,h,X)=u(kh,X) + h\times 64^m \, \left( \sum_{X \underset{m}\sim Y}  u(k\, h,Y)- u(k\,h,X)\right) \quad
\forall \, k\, \in\,\left \lbrace 0, \hdots, N-1 \right \rbrace,  \, \forall \,X\,\in \,V_m\setminus V_0
$$

\noindent We now fix~$m\,\in \,\N$, the number of contractions~$M$, and denote any~$X\,\in \,{\cal F}$ as~$X_{w,P_i}$, where~\mbox{$w\,\in\,\{1,\dots,N\}^m$} denotes a word of length $m$, and where $P_i$,~$0 \leq i \leq N_0-1$ belongs to~$V_0$. Let us also set:

$$n=\#\left \lbrace  w \, \big  | \quad \mid w\mid=m \right \rbrace \quad , \quad N_0=\#V_0$$

\noindent This enables one to introduce, for any integer~$k$ belonging to~$ \left \lbrace 0, \hdots, N-1 \right \rbrace$, the solution vector $U(k)$ as:

\[
U(k) =
\left(
\begin{matrix}
u(k\,h,X_{w_1 P_1})\\
u(k\,h,X_{w_1 P_2})\\
\vdots\\
u(k\,h,X_{w_2 P_1})\\
\vdots\\
u(k\,h,X_{w_n P_{N_0}})\\
\end{matrix}
\right)
\]

\noindent It satisfies the recurrence relation:

$$\forall\, k\,\in\, \left \lbrace 0, \hdots, N-1 \right \rbrace :\quad U(k+1) = A \, U(k)$$

\noindent where:

$$
A=I_{m}-h\times\tilde{\Delta}_{m}
$$

\noindent  and where~$I_{m}$ denotes the~$ {(\# V_m)}\times  {(\# V_m)}$ identity matrix, and $\tilde{\Delta}_{m}$ the~$ {(\# V_m)}\times  {(\# V_m)}$ normalized Laplacian matrix.\\

\vskip 1cm

\subsubsection{Theoretical study of the error, for H\"{o}lder continuous functions}

\noindent \emph{i}. \underline{General case}\\

In the spirit of the work of~R.~S.~Strichartz~\cite{Strichartz1999},~\cite{Strichartz2000}, it is interesting to consider the case of H\"{o}lder continuous functions. Why ? First, H\"{o}lder continuity implies continuity, which is a required condition for functions in the domain of the Laplacian (we refer to our work~\cite{RianeDavidM} for further details).\\
\noindent Second, a H\"{o}lder condition for such a function will result in fruitful estimates for its Laplacian, which is a limit of difference quotients.

\vskip 1cm

\noindent Let us thus consider a function~$u$ in the domain of the Laplacian, and a nonnegative real constant~$\alpha$ such that:

$$\forall\, (X,Y)\,\in \,{\cal F}^2,\forall\,t>0\, : \quad |u(t,X)-u(t,Y)| \leq C(t) \, |X-Y|^\alpha$$

\noindent where~$C$ denotes a positive function of the time variable~$t$.\\

\noindent Given a strictly positive integer~$m$, one has:

$$\Delta_m u(t,X) = \displaystyle\sum_{Y \in V_m,\,Y\underset{m}{\sim} X} \left (u(t,Y)-u(t,X)\right)  \quad \forall\,t>0,\, \forall\, X\,\in\, V_m\setminus V_0 $$

\noindent The Laplacian~$\Delta_{\mu} u$ is defined as the limit:

$$\Delta_{\mu} u(t,X)=\displaystyle\lim_{m\rightarrow+\infty} r^{-m} \left( \int_K{{\cal F}} \psi_{X_m}^{(m)} d\mu\right)^{-1}\, \Delta_m u(t, X_m)  \quad \forall\,t>0,\, \forall\, X\,\in\, K$$

\noindent where $\left (X_m\,\in\,V_m\setminus V_0 \right)_{m\in\N}$ is a sequence a points such that:
$$\displaystyle\lim_{m\rightarrow+\infty}X_m = X$$

\noindent and where~$r$ denotes the normalization ratio (we refer to~\cite{RianeDavidM}),~$\psi_{X_m}^{(m)}$ a harmonic spline function, and where:

$$\Delta_m u(t,X) = \displaystyle\sum_{Y \in V_m,\,Y\underset{m}{\sim} X} \left (u(t,Y)-u(t,X)\right)  \quad \forall\,t>0,\, \forall\, X\,\in\, V_m\setminus V_0 $$

\noindent Let us now introduce a strictly positive number~$\delta_{ij}=|P_i-P_j|$, for any~$P_i$ belonging to the set~\mbox{$V_0$}, and any~$P_j$ such that $P_j{\sim} P_i$. We set:
$$\delta_i=\max_{j}\delta_{ij}$$

\noindent We then introduce, for any integer~$i$ belonging to~$\left \lbrace   1, \hdots, M \right  \rbrace$, the contraction ratio of the similarity $f_i$, $R_i$, and set: $R=\displaystyle \max_{1 \leq i\leq M}R_i$.\\

\noindent One has then, for any~$X$ belonging to the set~\mbox{$ V_m\setminus V_0$}, any integer~$k$ belonging to~$\left \lbrace 0, \hdots, N-1 \right \rbrace$, and any strictly positive number~$h$:

$$\begin{array}{ccc} \left |r^{-m}\, \left( \displaystyle\int_K \psi_{X_m}^{(m)} d\mu\right)^{-1} \right | \, h \, \left |\Delta_m u(k\,h,X)\right | &\leq & \left | r^{-m} \, \left( \displaystyle\int_K \psi_{x_m}^{(m)} d\mu\right)^{-1} \right |\, h \, \displaystyle\sum_{Y \in V_m,\,Y\underset{m}{\sim} X} \left |u(k\,h,Y)-u(k\,h,X)\right| \\
&\leq & \left | r^{-m} \,\left( \displaystyle\int_K \psi_{X_m}^{(m)} d\mu\right)^{-1} \right | \,h\,C(k\,h) \,\displaystyle\sum_{Y \in V_m,\,Y\underset{m}{\sim} X}   |X-Y|^\alpha  \\
&\leq & \left | r^{-m}\,  \left( \displaystyle\int_K \psi_{X_m}^{(m)} d\mu\right)^{-1} \right | \,h\,C(k\,h) \,\displaystyle\sum_{m \, \mid \, Y \in V_m,\,Y\underset{m}{\sim} X}   \displaystyle \delta^{\alpha} \,R^{m\,\alpha}  \\
&\leq  & \left | r^{-m} \left( \int_K \psi_{X_m}^{(m)}\, d\mu\right)^{-1} \right | \,h\,C(k\,h) \,\displaystyle\sum_{p=0}^{+ \infty}   \displaystyle \delta^{\alpha}\, R^{p\, \alpha }  \\
&=  & \delta^{\alpha} \,  \displaystyle  \frac{ \left |r^{-m}\, \left( \displaystyle\int_K \psi_{X_m}^{(m)} d\mu\right)^{-1} \right | \,h\,C(k\,h)}{(1-R^{ \alpha })}  \\
\end{array}$$

\noindent We used the fact that, for $X\underset{m}{\sim} Y$, $X$ and $Y$ have addresses such that:

$$X=f_w(P_i) \quad , \quad Y=f_w(P_j)$$

\noindent for some $P_i$ and $P_j$ in $V_0$ and $w\in \{1,...,M\}^m$. May one set:

$$R(w) =R_{w_1}R_{w_2}...R_{w_m}$$

\noindent one gets:

$$\left |X-Y \right | =\left |fw(P_i)-f_w(P_j) \right |= R(w) \left |P_i-P_j \right |  \leq R^m \,\delta $$

\noindent The scheme~\mbox{$\left ({\cal S}_{\cal H}\right) $} allow us to write:

$$\begin{array}{ccccc}  \left |u((k+1)\,h,X)-u(k\,h,X) \right | &\leq &
 \delta^{\alpha} \,  \displaystyle  \frac{ \left |r^{-m}\, \left( \displaystyle\int_K \psi_{X_m}^{(m)} d\mu\right)^{-1} \right | \,h\,C(k\,h)}{(1-R^{ \alpha })}
\end{array}$$

\noindent One may note that a required condition for the convergence of the scheme is:
$$ \displaystyle \lim_{m\to + \infty,\, h \to 0^+} \left | r^{-m} \,\left( \displaystyle\int_K \psi_{X_m}^{(m)}\, d\mu\right)^{-1} \right | \,h\,C(k\,h)=0$$

\noindent In the case where~$C$ is a constant function, it reduces to:
$$ \displaystyle \lim_{m\to + \infty,\, h \to 0^+} \left | r^{-m} \,\left( \displaystyle\int_K \psi_{X_m}^{(m)}\, d\mu\right)^{-1} \right | \,h =0$$

\vskip 1cm
\newpage

\noindent \emph{ii.} \underline{The case of the Minkowski Curve~$\mathfrak MC$}\\

\noindent Given a strictly positive integer~$m$, due to:

$$\Delta_m u(t,X) = \displaystyle\sum_{Y \in V_m,\,Y\underset{m}{\sim} X} \left (u(t,Y)-u(t,X)\right)  \quad \forall\,t>0,\, \forall\, X\,\in\, V_m\setminus V_0 $$

\noindent one has then, for any~$X$ belonging to the set~\mbox{$ V_m\setminus V_0$}, any integer~$k$ belonging to~$\left \lbrace 0, \hdots, N-1 \right \rbrace$, and any strictly positive number~$h$:

$$\begin{array}{ccc}64^m\, h \, \left |\Delta_m u(k\,h,X)\right | &\leq &64^m\, h \, \displaystyle\sum_{Y \in V_m,\,Y\underset{m}{\sim} X} \left |u(k\,h,Y)-u(k\,h,X)\right| \\
&\leq &64^m\,h\,C(k\,h) \,\displaystyle\sum_{Y \in V_m,\,Y\underset{m}{\sim} X}   |X-Y|^\alpha  \\
&= &64^m\,h\,C(k\,h) \,\displaystyle\sum_{Y \in V_m,\,Y\underset{m}{\sim} X}   \displaystyle \frac{1 }{4^{m\,\alpha}}  \\
&\leq  &64^m\,h\,C(k\,h) \,\displaystyle\sum_{p=0}^{+ \infty}   \displaystyle \frac{1}{4^{p\, \alpha }}  \\
&=  &   \displaystyle  \frac{64^m\,h\,C(k\,h)}{1-4^{ - \alpha }}  \\
\end{array}$$

\noindent Using the scheme~\mbox{$\left ({\cal S}_{\cal H}\right) $}, one gets:

$$\begin{array}{ccccc}  \left |u((k+1)\,h,X)-u(k\,h,X) \right | &\leq &
 \displaystyle    \frac{64^m\,h\,C(k\,h)}{1-4^{ -\alpha  }}
\end{array}$$

\noindent One may note that a required condition for the convergence of the scheme is:
$$ \displaystyle \lim_{m\to + \infty,\, h \to 0^+} 64^m\,h\,C(k\,h)=0$$

\noindent In the case where~$C$ is a constant function, it reduces to:
$$ \displaystyle \lim_{m\to + \infty,\, h \to 0^+} 64^m\,h =0$$

\vskip 1cm

\newpage

\subsubsection{Numerical results}
\noindent In the sequel, we present the numerical results for~$m=3$, $T=10$ and $N=10^7$. One may note that, as expected, the solution shows abnormal behavior until the ration $h\times 64^m$ goes towards zero.

\begin{center}
\includegraphics[scale=1.5]{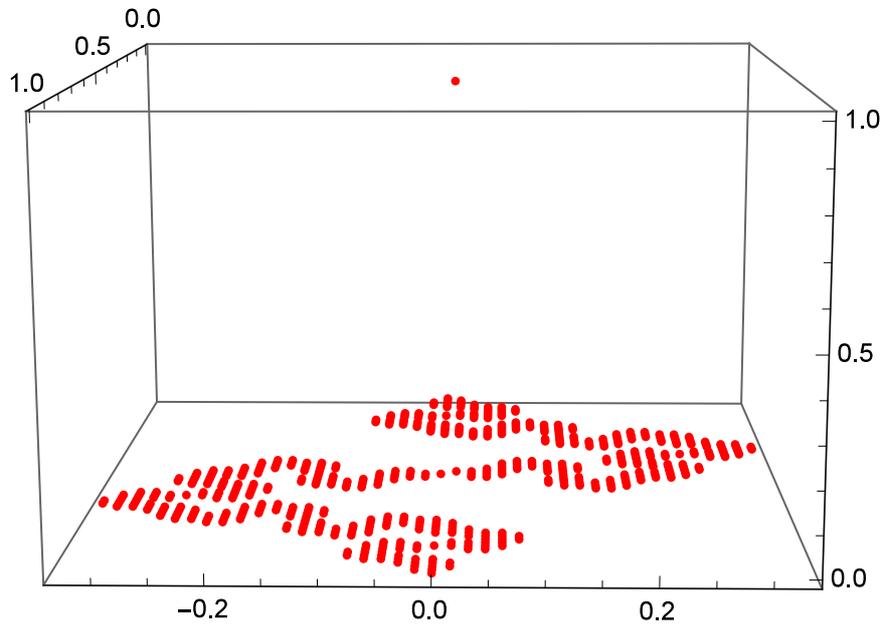}
\captionof{figure}{The graph of the approached solution of the heat equation for $t=0$.}
\label{fig1}
\end{center}

\begin{center}
\includegraphics[scale=1.5]{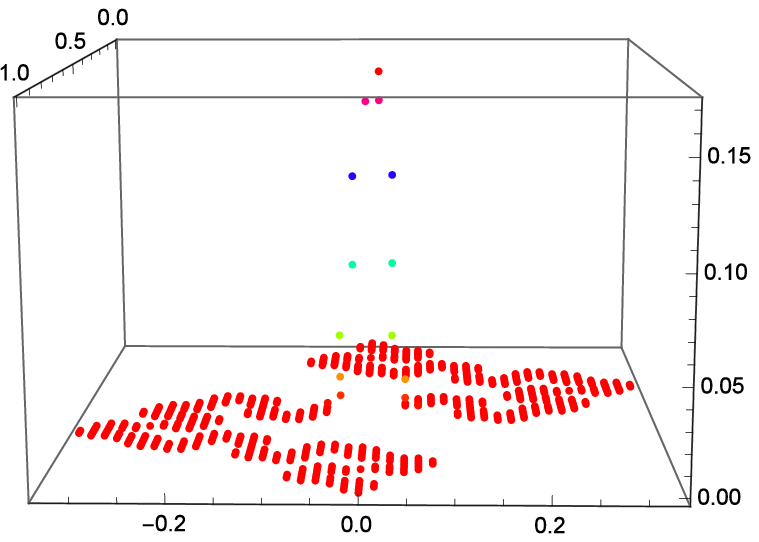}
\captionof{figure}{The graph of the approached solution of the heat equation for $t=10$.}
\label{fig1}
\end{center}

\begin{center}
\includegraphics[scale=1.5]{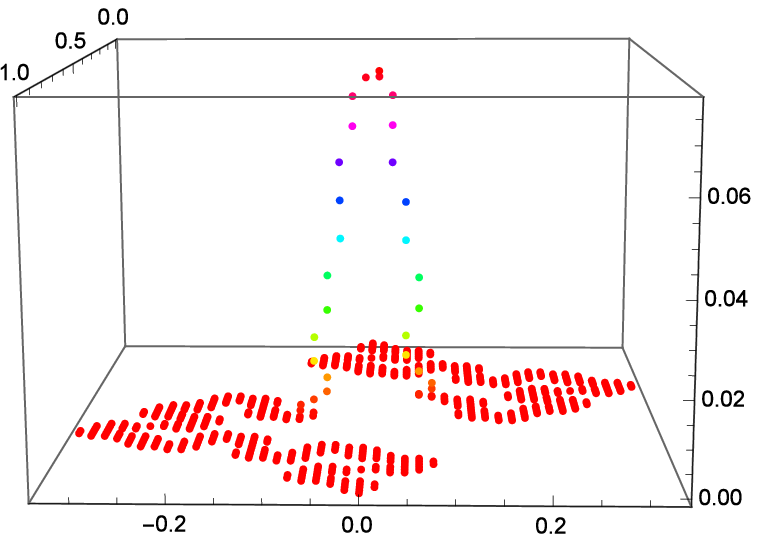}
\captionof{figure}{The graph of the approached solution of the heat equation for $t=50$.}
\label{fig1}
\end{center}

\begin{center}
\includegraphics[scale=1.5]{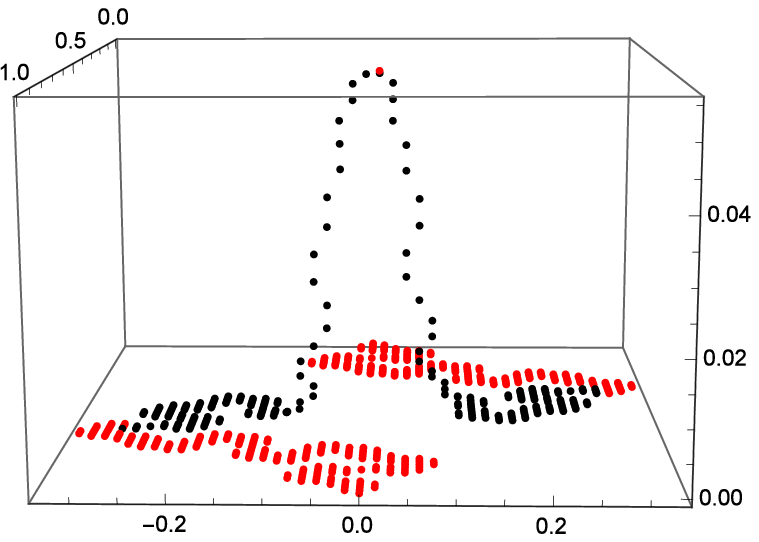}
\captionof{figure}{The graph of the approached solution of the heat equation for $t=100$.}
\label{fig1}
\end{center}

\begin{center}
\includegraphics[scale=1.5]{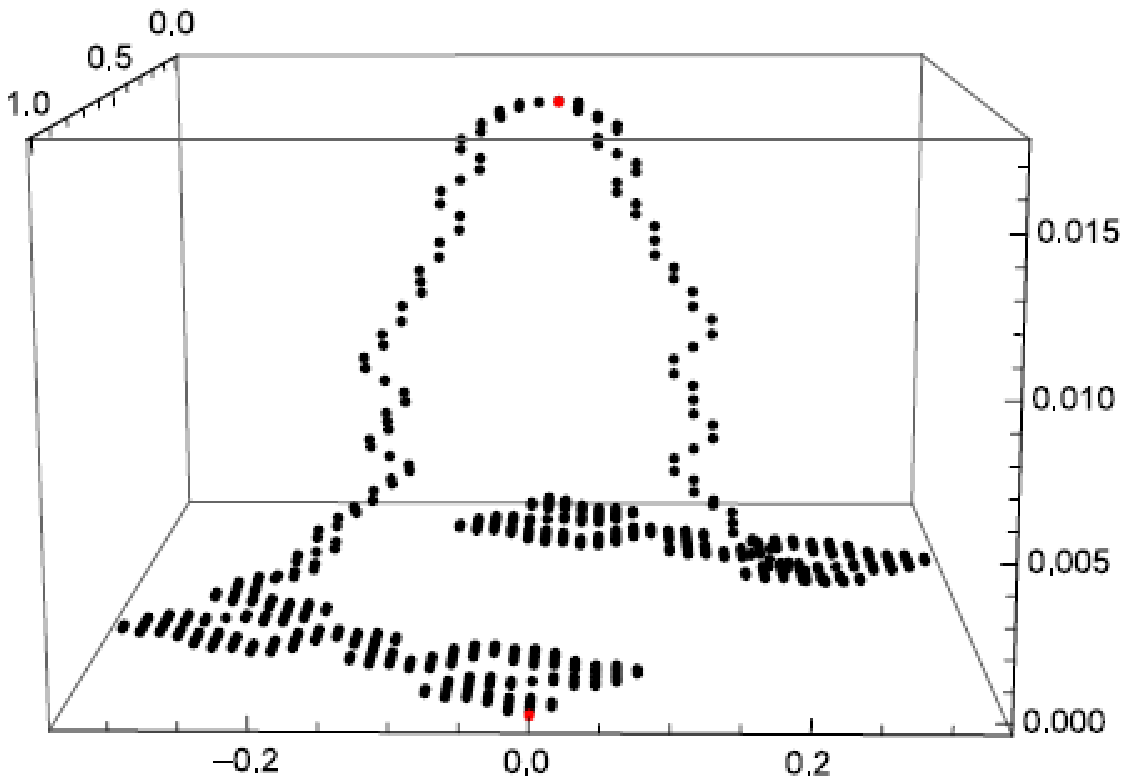}
\captionof{figure}{The graph of the approached solution of the heat equation for $t=1000$.}
\label{fig1}
\end{center}

\newpage
\subsection{The wave equation}

\subsubsection{Formulation of the problem}

One may proceed in a similar way as in the above. First, we consider a solution~$u$ of the problem:

$$ \left \lbrace  \begin{array}{ccccc}
\displaystyle \frac{\partial u}{\partial t}(t,x)-\Delta u(t,x)&=&0 & \forall (t,x)\,\in \,\left]0,T\right[  \times {\cal F}\\
u(t,x)&=&0 & \forall \, ( x,t) \,\in \,\partial {\cal F} \times \left[0,T\right[\\
u(0,x)&=&g(x)  & \forall x\,\in \,{\cal F}
\end{array}\right.$$

\noindent In order to define a numerical scheme, one may use a second order central difference in time. The Laplacian is approximated by means of the graph Laplacians~$\Delta_m \,u$, defined on the sequence of graphs~$\left ( {\mathfrak MC}_m\right)_{m\in\N}$. \\

\noindent To this purpose, we fix a strictly positive integer~$N $, and set:

$$ \displaystyle{h =\frac{T}{N}}$$

\noindent One has, for any integer~$k$ belonging to~$\left \lbrace 0, \hdots, N-1 \right \rbrace$:

$$   \forall \,X\,\in \,{\cal F} \,:   \quad
\displaystyle\frac{\partial u}{\partial t}(kh,x)=\displaystyle\frac{1}{h}\,\left( u((k+1)\,h,X)-u(kh,x)\right)+{\cal O}(h)
$$

\noindent and:
$$   \forall \,X\,\in \,{\cal F} \,:   \quad
\displaystyle\frac{\partial^2 u}{\partial t^2}(k\,h,X)=\displaystyle\frac{ u((k+1)\,h,X)-2u(k\,h,X)+u((k-1)\,h,X)}{h^2}+O(h^2)
$$

\noindent According to~\cite{RianeDavidM}, the Laplacian on the Minkowski curve is given by:

 $$ \forall X\,\in \, {\cal F} \,:   \quad
\Delta u(t,X)= \lim_{m\rightarrow +\infty} 64^m \, \left( \sum_{X \underset{m}\sim Y}  u(t,Y)- u(t,X)\right)
$$

\noindent This enables one to approximate the Laplacian, at a~$m^{th}$ order,~$m \,\in\,\N^\star$, using the graph normalized Laplacian as follows:

$$ \forall \, k\, \left \lbrace 0, \hdots, N-1 \right \rbrace,\,\forall \,X\,\in \, {\cal F}   \,: \quad
\Delta u(t,X)\approx 64^m \,\left( \sum_{X \underset{m}\sim y}  u(k\,h,Y)- u(kh,X)\right)
$$

\noindent One thus obtains the following scheme, for any integer~$k$ belonging to~$\left \lbrace 0, \hdots, N-1 \right \rbrace$, any point~$P_j$ of~$V_0$,~$0\leq j \leq N_0-1$, and any~$X$ in the set~\mbox{$V_m\setminus V_0$}:

$$\left ({\cal S}_{\cal W}\right) \quad\left \lbrace \begin{array}{cccc}
\displaystyle\frac{u((k+1)\,h,X)-2u(k\,h,X)+u((k-1)\,h,X)}{h^2}&=&64^m \,\left( \displaystyle\sum_{X \underset{m}\sim Y}  u(k\,h,Y)- u(k\,h,X)\right) &  \\
u(k\,h,P_j)&=&0 &  \\
u(X,0)&=&g(X) &
\end{array}\right.
$$
\normalsize

\noindent Let us define the approximate equation as:
$$
u((k+1)\,h,X)=2\, u(k\,h,X)-u((k-1\,)h,X)+ h^2\times 64^m \left( \sum_{X \underset{m}\sim Y}  u(k\,h,Y)- u(k\,h,X)\right)
$$

\noindent We now fix~$m\,\in \,\N$, the number of contractions~$M$  and denote any~$X\,\in \,K$ as~$X_{w,P_i}$, where~\mbox{$w\,\in\,\{1,\dots,N\}^m$} denotes a word of length $m$, and where $P_i$,~$0 \leq i \leq N_0-1$ belongs to~$V_0$. Let us also set:

$$n=\#\left \lbrace  w \, \big  | \quad \mid w\mid=m \right \rbrace \quad , \quad N_0=\#V_0$$

\noindent This enables one to introduce, for any integer~$k$ belonging to~$ \left \lbrace 0, \hdots, N-1 \right \rbrace$, the solution vector $U(k)$ as:

\[
U(k) =
\left(
\begin{matrix}
u(k\,h,X_{w_1 P_1})\\
u(k\,h,X_{w_1 P_2})\\
\vdots \\
u(k\,h,X_{w_2 P_1})\\
\vdots\\
u(k\,h,X_{w_n P_{N_0}})\\
\end{matrix}
\right)
\]

\noindent One has:

$$U(k+1) = A\, U(k)-U(k-1)$$

\noindent where:

$$
A=2\,I{N_0}-h^2 \times\tilde{\Delta}_{m}
$$

\noindent  and where~$I_{m}$ denotes the~${(\#V_m)}\times  {(\#V_m)}$ identity matrix, and $\tilde{\Delta}_{m}$ the~${(\#V_m)}\times  {(\#V_m)}$ normalized Laplacian matrix.\\

\subsubsection{Theoretical study of the error, for H\"{o}lder continuous functions}

A previously, let us thus consider a function~$u$ in the domain of the Laplacian, and a nonnegative real constant~$\alpha$ such that:

$$\forall\, (X,Y)\,\in {\mathfrak {MC} }^2,\forall\,t>0\, : \quad |u(t,X)-u(t,Y)| \leq C(t) \, |X-Y|^\alpha$$

\noindent where~$C$ denotes a positive function of the time variable~$t$.\\

\noindent Given a strictly positive integer~$m$, due to:

$$\Delta_m u(t,X) = \displaystyle\sum_{Y \in V_m,\,Y\underset{m}{\sim} X} \left (u(t,Y)-u(t,X)\right)  \quad \forall\,t>0,\, \forall\, X\,\in\, V_m\setminus V_0 $$

\noindent one has then, for any~$X$ belonging to the set~\mbox{$ V_m\setminus V_0$}, any integer~$k$ belonging to~$\left \lbrace 0, \hdots, N-1 \right \rbrace$, and any strictly positive number~$h$:

$$\begin{array}{ccc}64^m\,h^2\, \left |\Delta_m u(k\,h,X)\right | &\leq   &64^m\,h^2\,C(k\,h) \,\displaystyle  \frac{1}{1-4^{ - \alpha }}  \\
\end{array}$$

\noindent Using the scheme~\mbox{$\left ({\cal S}_{\cal W}\right) $}, one gets:

$$\begin{array}{ccccc}  \left |u((k+1)\,h,X)-2u(k\,h,X)+u((k-1)\,h,X)) \right | &\leq &
 \displaystyle    \frac{64^m\,h^2\,C(k\,h)}{1-4^{ -\alpha  }}
\end{array}$$

\noindent One may note that a required condition for the convergence of the scheme is:
$$ \displaystyle \lim_{m\to + \infty,\, h \to 0^+} 64^m\,h^2\,C(k\,h)=0$$

\noindent In the case where~$C$ is a constant function, it reduces to:
$$ \displaystyle \lim_{m\to + \infty,\, h \to 0^+} 64^m\,h^2 =0$$

\vskip 1cm

\subsubsection{Numerical results}

\noindent In the sequel, we present the numerical results for~$m=2$, $T=10$ and $N=10^3$. As it was the case for the heat equation, the solution shows abnormal behavior until the ration $h^2\times64^m$ goes towards zero.

\begin{center}
\includegraphics[scale=1.5]{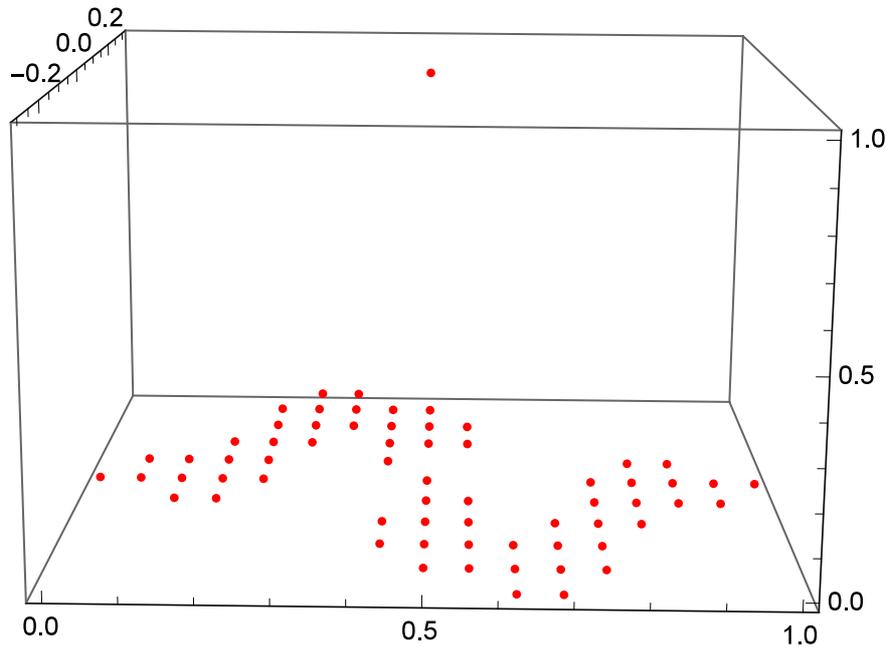}
\captionof{figure}{The graph of the approached solution of the wave equation for $t=0$.}
\label{fig1}
\end{center}

\begin{center}
\includegraphics[scale=1.5]{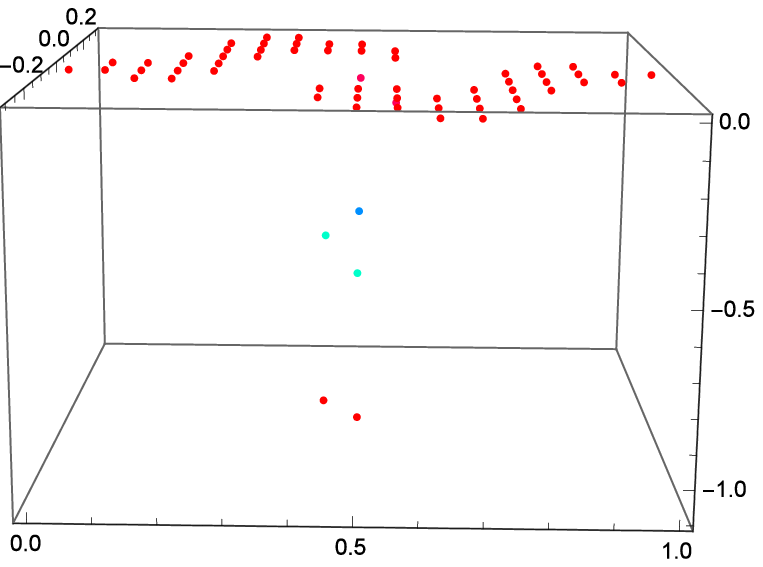}
\captionof{figure}{The graph of the approached solution of the wave equation for $t=5$.}
\label{fig1}
\end{center}

\begin{center}
\includegraphics[scale=1.5]{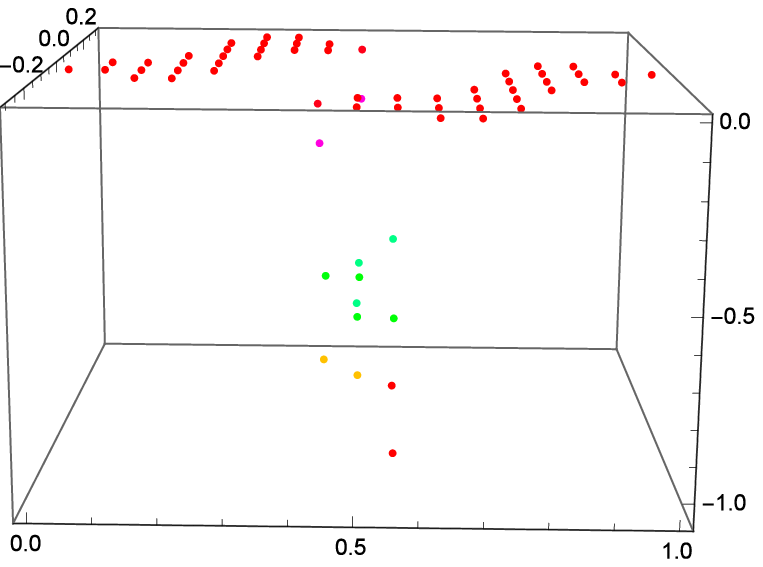}
\captionof{figure}{The graph of the approached solution of the wave equation for $t=10$.}
\label{fig1}
\end{center}

\begin{center}
\includegraphics[scale=1.5]{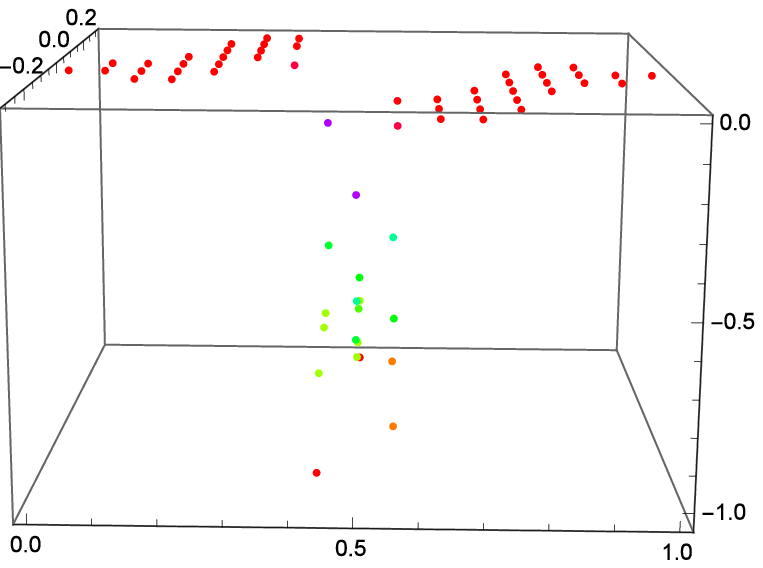}
\captionof{figure}{The graph of the approached solution of the wave equation for $t=15$.}
\label{fig1}
\end{center}

\begin{center}
\includegraphics[scale=1.5]{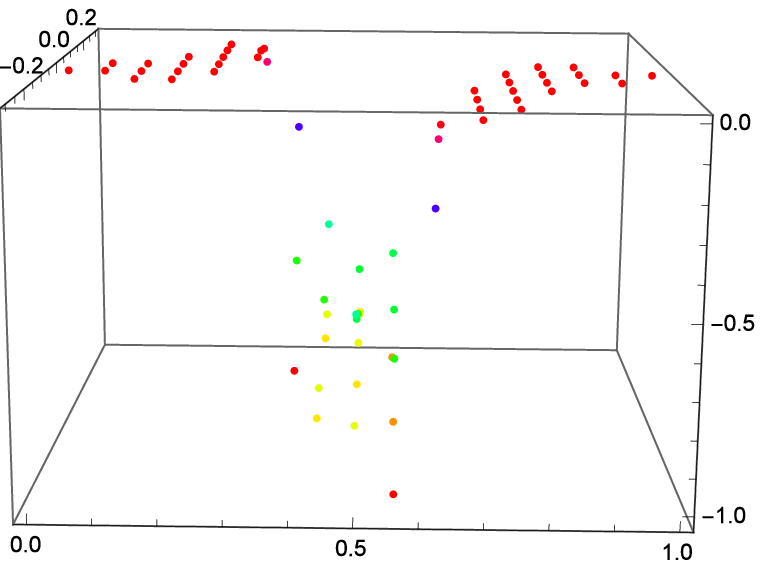}
\captionof{figure}{The graph of the approached solution of the wave equation for $t=20$.}
\label{fig1}
\end{center}

\subsection{Error computing}

\noindent Let us recall that, given a first order forward difference in time, a strictly positive integer~$N $, and a strictly positive real number~$h =\displaystyle \frac{T}{N}$, one has:

$$\forall \, X \in\,{\mathfrak MC} \, : \quad
\displaystyle \frac{\partial u}{\partial t}(k\,h,X)-\frac{1}{h}\left( u((k+1)\,h,X)-u(k\,h,X)\right)={ \cal O}(h)
$$

\noindent and:

$$\forall \, X \in\,{\cal F} \, : \quad
\displaystyle \frac{\partial^2 u}{\partial t^2}(k\,h,X)-\displaystyle \frac{1}{h^2}\left( u((k+1)\,h,X)-2\,u(k\,h,X)+u((k-1)\,h,X)\right)={\cal O}(h^2)
$$

\noindent If one aims at computing the error related to the use of discrete Laplacian, one may consider the Dirichlet problem:

$$\left \lbrace \begin{array}{ccc}
-\Delta u + q\, u &= &f\\
u_{\mid \partial {\cal F}}&=&0
\end{array}\right.$$

\noindent We choose $u$ to be harmonic, and we fix $q=2$ to obtain $f=2\,u$. We calculate next the solution corresponding to the finite difference scheme for this Dirichlet problem.\\

\noindent As previously, the Laplacian is approximated by means of the graph Laplacians~$\Delta_m \,u$, defined on the sequence of graphs~$\left ( {\mathfrak MC}_m\right)_{m\in\N}$:

$$ \forall \, k\, \left \lbrace 0, \hdots, N-1 \right \rbrace,\,\forall \,(t,X)\,\in \, [0,T] \times  {\mathfrak MC} \,: \quad
\Delta u(t,X)= 64^m \,\left( \displaystyle \sum_{X \underset{m}\sim Y}  u(k\,h,Y)- u(k\,h,X)\right) + E_m
$$

\noindent where $E_m$ denotes an error term, that may be obtained numerically. One gets:

$$\left \lbrace \begin{array}{cccc}
-64^m \,\left( \displaystyle \sum_{X \underset{m}\sim Y}  u(Y)- u(X)\right)+q u(X)&=&f(X) &  \forall \,X\, \in \, V_m \setminus V_0\\
u(P_j)&=&0 & \forall \, P_j \, \in \, V_0\\

\end{array}\right.$$

\noindent One defines the approximate equation as:

\[
-64^m \, \left( \displaystyle\sum_{X \underset{m}\sim Y}  u(Y)- u(C)\right)+q \, u(X)=f(X)
\]

\noindent We now fix $m\,\in \,\N$, and denote any~$C\,\in \,{\cal F}$ as~$X_{w,P_i}$, where~$w\in\{1,\dots,N\}^m$ denotes a word of length $m$, and where $P_i$,~$0 \leq i \leq N_0-1$ belongs to~$V_0$. Let us also set:

$$n=\#\left \lbrace  w \, \big  | \quad \mid w\mid=m \right \rbrace \quad , \quad N_0=\#V_0$$

\noindent We then introduce the vector $U$:

\[
U =
\left(
\begin{matrix}
u(X_{w_1 P_1})\\
u(X_{w_1 P_2})\\
\vdots\\
u(X_{w_2 P_1})\\
\vdots\\
u(X_{w_n P_{N_0}})\\
\end{matrix}
\right)
\]

\noindent One has:

$$ A\, U=F$$

\noindent where:

$$
A= \tilde{\Delta}_{m}+q\, I_{m}
$$

\noindent and:

$$
F =
\left(
\begin{matrix}
f(X_{w_1 P_1})\\
f(X_{w_1 P_2})\\
\vdots\\
f(X_{w_2 P_1})\\
\vdots\\
f(X_{w_n P_{N_0}})\\
\end{matrix}
\right)
$$

\noindent where $I_{m}$ denotes the~${(\#V_m)}\times  {(\#V_m)}$ identity matrix, and $\tilde{\Delta}_{m}$ the~${(\#V_m)}\times  {(\#V_m)}$ normalized Laplacian matrix.

\noindent One may begin by evaluating the error of the Laplacian approximation on the Minkowski curve, we use the Dirichlet problem :

$$ \left \lbrace \begin{array}{ccc}
-\Delta u + q \,u &= &f\\
u_{\mid \partial {\mathfrak MC}}&=&0
\end{array}\right.$$

\noindent for an harmonic function $u$, the problem has a solution. The error is then:
$$E_m =\parallel u-u_m\parallel_{\infty}$$

%\noindent The following figures show the variations of the error as a function of the graph order $m$:

%\begin{center}
%\includegraphics[scale=1.5]{MC-error.eps}
%\captionof{figure}{The scheme error as a function of~$m$.}
%\label{fig1}
%\end{center}

\section{Conclusion}

\hskip 0.5cm Refined numerical methods that fit problems of everydaylife like diffusion, or transmission of waves by fractal micro-antennas, seem unavoidable. As for now, current computations do not take into account advances in fractal analysis.\\
Our finite-difference scheme is destined to help improve those calculations. Lots remain to be done: in a near future, specific tests are to be made with existing Minkowski fracatal antennas.

\section*{Annex: Algorithm and Mathematica program}

%See the file FDM-MI.nb.\\
In the sequel, we explain how one may implement our algorithm. We choose to build a specific Mathematica program, that enables one to compute the numerical solution of the heat equation. The program break downs into five parts, which arise naturally from the above theoretical study.\\

\subsection*{Annex 1. Initialization}

One here defines the initial points~$P_0$ and $P_1$ that constitutes the set~$V_0$, as well as the similarities $f_i$,~\mbox{$ 0 \leq i \leq N-1$}.\\

\subsection*{Annex 2. Harmonic functions}

One here defines a harmonic function "HadresseH" on~$\mathfrak{MC}$ ; given the value of a harmonic function $f$ such that:
$$f(P_0) = a \,\in\,\R \quad , \quad  f (P_1) = b \,\in\,\R$$

\noindent One thus gets the value of $f$ at a point, the address of which~$w$ belongs to~$\left \lbrace  1, ..., 8 \right \rbrace ^m$,~$m \,\in\,\N^\star$, for an initial point~\mbox{ $p\, \in \,V_0$.} For instance:

$$\text{HadresseH}[0, 1, {1, 2, 1}, p0]=\displaystyle \frac{1}{64}$$

\subsection*{Annex 3. Construction of the Laplacian matrix}

\noindent We build a function "sparseMat" enables one to construct a~$ n \times n$ tridiagonal matrix,~$n \,\in\,\N^\star$, with respective values $x$,~$y$ and $z$ on the tridiagonal ($x$ on the diagonal).\\

\noindent We have also built a function "MatL" gives the Laplacian matrix of order~$n$, i.e., associated with the~$n^{th}$ order approximate graph. The matrix is also a tridiagonal one, since each point has a unique predecessor and a single successor with the exception of the first and last point ($ V_0 $ is the boundary). For instance:

$$\text{MatL[1] }=\left(
\begin{array}{ccccccccc}
 1 & -1 & 0 & 0 & 0 & 0 & 0 & 0 & 0 \\
 -1 & 2 & -1 & 0 & 0 & 0 & 0 & 0 & 0 \\
 0 & -1 & 2 & -1 & 0 & 0 & 0 & 0 & 0 \\
 0 & 0 & -1 & 2 & -1 & 0 & 0 & 0 & 0 \\
 0 & 0 & 0 & -1 & 2 & -1 & 0 & 0 & 0 \\
 0 & 0 & 0 & 0 & -1 & 2 & -1 & 0 & 0 \\
 0 & 0 & 0 & 0 & 0 & -1 & 2 & -1 & 0 \\
 0 & 0 & 0 & 0 & 0 & 0 & -1 & 2 & -1 \\
 0 & 0 & 0 & 0 & 0 & 0 & 0 & -1 & 1 \\
\end{array}
\right)$$

\subsection*{Annex 4. Construction of the initial vector}

\noindent We build a function called "Genadresse", that defines any required address, the length of which take the value~$n\,\in\,\N^\star$.\\

\noindent We have also introduced a function "Fadress", which enables one to calculate the coordinates of a point, given its address~$w$, starting from an initial point~$p$.\\

\noindent To define the initial vector, i.e. the one that yields all the points of the approximate graph~${\mathfrak MC}_m$ of order $ m\,\in\,\N^\star $, one uses the function "InitV". We have used the fact that any point of $ V_m \setminus V_0 $ owns two addresses:
$$f_{w_1 ... w_m} (P_0) = f_{w_1 ... (w_m+1)} (P_1)$$

\noindent This enables one to obtain the points~$V_m$ by applying all the similarities of order~$m$ to the point $P_1$, and, then, add the missing point~$ P_0 $. For instance:

$$\text{InitV[1]}=\left \lbrace \{0, 0\}, \left \lbrace\displaystyle \frac{1}{4}, 0\right \rbrace ,  \left \lbrace\displaystyle \frac{1}{4}, \displaystyle \frac{1}{4}\right \rbrace,
 \left \lbrace\{\displaystyle \frac{1}{2}, \displaystyle \frac{1}{4}\right \rbrace,\left \lbrace\displaystyle \frac{1}{2}, 0\right \rbrace,
\left \lbrace\displaystyle \frac{1}{2}, -\displaystyle \frac{1}{4}\right \rbrace, \left \lbrace\displaystyle \frac{3}{4}, -\displaystyle \frac{1}{4}\right \rbrace,
\left \lbrace\displaystyle \frac{3}{4}, 0\right \rbrace, \{1, 0\}\right \rbrace $$

\subsection*{Annex 5. Solution of the Heat equation}

\noindent One has to define an initial function~$gH$, which corresponds to the value of the searched solution~$u$, at the initial time~\mbox{$t = 0$}. For the sake of simplicity, we choose the identically null function~$gH=0$.\\

\noindent The function "AAH [N, T, k, n]" yields the matrix $I-h\, \tilde{\Delta_m}$ for a solution with horizon ??? $T$ ; $N$ denotes the number of steps, $k$, the order of the approximate graph ${\mathfrak{MC}}_k$, and $n$ the number of points. For instance:

$$\text{AAH[10, 10, 1, 9]}= \left(
\begin{array}{ccccccccc}
 -63 & 64 & 0 & 0 & 0 & 0 & 0 & 0 & 0 \\
 64 & -127 & 64 & 0 & 0 & 0 & 0 & 0 & 0 \\
 0 & 64 & -127 & 64 & 0 & 0 & 0 & 0 & 0 \\
 0 & 0 & 64 & -127 & 64 & 0 & 0 & 0 & 0 \\
 0 & 0 & 0 & 64 & -127 & 64 & 0 & 0 & 0 \\
 0 & 0 & 0 & 0 & 64 & -127 & 64 & 0 & 0 \\
 0 & 0 & 0 & 0 & 0 & 64 & -127 & 64 & 0 \\
 0 & 0 & 0 & 0 & 0 & 0 & 64 & -127 & 64 \\
 0 & 0 & 0 & 0 & 0 & 0 & 0 & 64 & -63 \\
\end{array}
\right)$$

\noindent One then defines a function "Per[k, j, m]" which, given a graph of order $ k $ and a number of points $m$, associates the value $1$ to the point of order $j$ in the "initial vector".\\

\noindent Another function, "UH[t, k, T, N, g, P, n, m]" recursively gives, at time~$ t >0 $, for a graph order $k$, a horizon~$T$, and a step number~$N$, the value reached by an initial function $g + P$.
% (one must bear in mind that the initial function is the sum of the null function and the perturbation), a perturbation order $m$ and a number of points $n$, the solution ??? %sens ???. \\

\noindent A function "DFHA[t, n, T, N, g, P, j]" yields the coordinate triplet $(x, y)$ $+$ solution.

\bibliographystyle{alpha}
%\bibliography{BibliographieClaire}

\end{document}